\documentclass[12pt,fleqn]{article}
\usepackage{hyperref}
\usepackage{amsfonts}
\usepackage{amssymb}
\usepackage{graphicx}
\newcommand{\ol}{\setlength{\itemsep}{0pt.}\begin{enumerate}}
\newcommand{\eol}{\end{enumerate}\setlength{\itemsep}{-\parsep}}
\newcommand{\ignore}[1]{}
\newcommand{\quotes}[1]{``#1''}

\newtheorem{THEOREM}{Theorem}[section]
\newenvironment{theorem}{\begin{THEOREM} \hspace{-.85em} {\bf :}
}%
                        {\end{THEOREM}}
\newtheorem{LEMMA}[THEOREM]{Lemma}
\newenvironment{lemma}{\begin{LEMMA} \hspace{-.85em} {\bf :} }%
                      {\end{LEMMA}}
\newtheorem{COROLLARY}[THEOREM]{Corollary}
\newenvironment{corollary}{\begin{COROLLARY} \hspace{-.85em} {\bf
:} }%
                          {\end{COROLLARY}}
\newtheorem{PROPOSITION}[THEOREM]{Proposition}
\newenvironment{proposition}{\begin{PROPOSITION} \hspace{-.85em}
{\bf :} }%
                            {\end{PROPOSITION}}
\newtheorem{DEFINITION}[THEOREM]{Definition}
\newenvironment{definition}{\begin{DEFINITION} \hspace{-.85em} {\bf
:} \rm}%
                            {\end{DEFINITION}}
\newtheorem{EXAMPLE}[THEOREM]{Example}
\newenvironment{example}{\begin{EXAMPLE} \hspace{-.85em} {\bf :}
\rm}%
                            {\end{EXAMPLE}}
\newtheorem{CONJECTURE}[THEOREM]{Conjecture}
\newenvironment{conjecture}{\begin{CONJECTURE} \hspace{-.85em}
{\bf :} \rm}%
                            {\end{CONJECTURE}}
\newtheorem{MAINCONJECTURE}[THEOREM]{Main Conjecture}
\newenvironment{mainconjecture}{\begin{MAINCONJECTURE} \hspace{-.85em}
{\bf :} \rm}%
                            {\end{MAINCONJECTURE}}
\newtheorem{PROBLEM}[THEOREM]{Problem}
\newenvironment{problem}{\begin{PROBLEM} \hspace{-.85em} {\bf :}
\rm}%
                            {\end{PROBLEM}}
\newtheorem{QUESTION}[THEOREM]{Question}
\newenvironment{question}{\begin{QUESTION} \hspace{-.85em} {\bf :}
\rm}%
                            {\end{QUESTION}}
\newtheorem{REMARK}[THEOREM]{Remark}
\newenvironment{remark}{\begin{REMARK} \hspace{-.85em} {\bf :}
\rm}%
                            {\end{REMARK}}
\newtheorem{CLAIM}[THEOREM]{Claim}
\newenvironment{claim}{\begin{CLAIM} \hspace{-.85em} {\bf :}
\rm}%
                            {\end{CLAIM}}
\newcommand{\thm}{\begin{theorem}}
\newcommand{\lem}{\begin{lemma}}
\newcommand{\pro}{\begin{proposition}}
\newcommand{\dfn}{\begin{definition}}
\newcommand{\rem}{\begin{remark}}
\newcommand{\xam}{\begin{example}}
\newcommand{\cnj}{\begin{conjecture}}
\newcommand{\mcnj}{\begin{mainconjecture}}
\newcommand{\prb}{\begin{problem}}
\newcommand{\que}{\begin{question}}
\newcommand{\cor}{\begin{corollary}}
\newcommand{\clm}{\begin{claim}}

\newcommand{\ethm}{\end{theorem}}
\newcommand{\elem}{\end{lemma}}
\newcommand{\epro}{\end{proposition}}
\newcommand{\edfn}{\end{definition}}
\newcommand{\erem}{\bbox\end{remark}}
\newcommand{\exam}{\bbox\end{example}}
\newcommand{\ecnj}{\bbox\end{conjecture}}
\newcommand{\emcnj}{\bbox\end{mainconjecture}}
\newcommand{\eprb}{\bbox\end{problem}}
\newcommand{\eque}{\bbox\end{question}}
\newcommand{\ecor}{\end{corollary}}
\newcommand{\eclm}{\end{claim}}

\newcommand{\beqn}{\begin{equation}}
\newcommand{\eeqn}{\end{equation}}

\newcommand{\bbox}{\begin{flushright} $\Box $ \end{flushright}}
\newcommand{\qed}{\bbox}

\def\eps{\epsilon}

\def \1{\mathbf 1}

\newcommand{\CC}{{\cal C}}

\def\<{\left<}
\def\>{\right>}
\def \({\left(}
\def \){\right)}

\def \8{\infty}

\newcommand{\remove}[1]{}

\title{The success probability in Levine's hat problem, and independent sets in graphs }
\author{Noga Alon \and Ehud Friedgut \and Gil Kalai\and Guy Kindler}
\begin{document}
\maketitle
\date{}
\begin{abstract}
Lionel Levine's hat challenge has $t$ players, each with a (very large, or infinite) stack of hats on their head, each hat independently colored at random black or white. 
The players are allowed to coordinate before the random colors are chosen, but not after.
Each player sees all hats except for those on her own head.
They then proceed to simultaneously try and each pick a black hat from their respective stacks.
They are proclaimed successful only if they are all correct.
Levine's conjecture is that the success probability tends to zero when the number of players grows. 
We prove that this success probability is strictly decreasing in the number of players, and present some connections to problems in graph theory: relating the size of the largest independent set in a graph and in a random induced subgraph of it, and bounding the size of a set of vertices intersecting every maximum-size independent set in a graph.
\end{abstract}

\section{Introduction}
The following question proposed by Lionel Levine, arose in the context of his work with Friedrich \cite{FL13}. It gained considerable popularity after being presented in 2011 Tanya Khovanova's blog \cite{Khova}.
Consider $t$ players, each with a stack of $n$ hats on her head, where the hats are chosen independently at random to be either black or white with probability 1/2. Each player sees the hats of every other player, but not her own. Then, simultaneously, all players pick a hat from their respective stacks. The collective of players wins if every single player points to a black hat, else, if even a single player errs, the collective fails. Let $p(t,n)$ be the maximal success probability over all possible strategies that the players can apply. Let $p(t)$ be the limit of $p(t,n)$ as $n$ tends to infinity. The challenge set by Levine was to prove the following conjecture.
\begin{conjecture}\label{Main}
$p(t)$ tends  to 0 as $t$ grows.
\end{conjecture}
In this paper we prove
\begin{theorem}\label{thm:main}
 $p(t+1) < p(t)$ for all $t\geq1$.
 \end{theorem}
 While preparing this paper we were sent a draft of a comprehensive hat-related paper by Buhler, Freiling, Graham, Kariv, Roche, Tiefenbruck, Van Alten and Yeroshkin \cite{BFGKRTVY}, where Theorem \ref{thm:main} is also proven, along with other interesting results and bounds. We refer to their paper as an excellent source for background on the state of the art for this problem. The most prominent landmarks mentioned in their paper are 
$$
0.35 \leq p(2) \leq 0.361607
$$
and
$$
p(t) = \Omega (1/\log(t)),
$$
where the bounds on $p(2)$ are due to them, and the $1/\log(t)$ probably due to Peter Winkler.
The fact that $p(2) \leq 3/8$ is well known folklore in the hatter community. Alon and Tardos have an approach to improve it \cite{AT}, but not to as tight a bound as 0.362.

 In this paper, we also present some generalizations of Levine's conjecture, relating it to questions regarding independent sets in Hamming products of graphs, and independent sets in random induced subgraphs. In this context we prove that for every graph $G$ on $n$ vertices with independence
number $(1/4+\eps)n$, the average independence number of an induced
subgraph of $G$ on a uniform random subset of the vertices is at
most $(1/4+\eps-\Omega(\eps^2)) n$. For regular graphs we improve the range of these bounds to graphs with  independent sets of size slightly larger than $n/8$. These are theorems \ref{t16} and \ref{p17}
 
 The proof of our main theorem leads naturally to the question of bounding the size of a set which intersects all large maximal (with respect to containment) independent sets, or all maximum-size  independent sets in a graph, and is related to a conjecture regarding this by Bollob\'as, Erd\H{o}s and Tuza, Conjecture \ref{cbet}. 
 We give a construction yielding a bound related to this conjecture - an infinite family of graphs $G_n$, where
$G_n$ has $n$ vertices, independence number at least $n/4$, and no set of less than
$\sqrt{n}/2$ vertices intersects all its maximum independent sets.

\section{General setting, strategies and winning sets}
Let us start by defining a general setting that includes the hats game as a special case.
Let $B$ be a fixed ground set, and let $\cal{W}$ a family of subsets of $B$, that we will call winning sets (In the case of the hats game, $B=\{0,1\}^n$, and the winning sets are the half-cubes of the form $\{x: x_i=1\}$).  The corresponding game there are $t \geq 1$ players, each is assigned a point at random from $B$, and each player sees the points the other players were assigned, but not her own point (which is \quotes{on her forehead}). Then, simultaneously, each player chooses a winning set, and the collective of players succeeds if every player named a set containing her point.

 Let us now define what a strategy is for $B^t$ (\quotes{the game for $t$ players}), and what a winning set is for such a strategy. For $t=1$ a winning set is any of the sets in  $\cal{W}$, and a strategy is a choice of one winning set, namely, a function $f: \{\emptyset\} \rightarrow$  $\cal{W}$, so the winning set for the strategy $f$ is $f(\emptyset)$.)
 For the $t$-player game a strategy is a $t$-tuple of functions, 
 $(f_1,\ldots,f_t)$, where $f_i : B^{t-1} \rightarrow$ $\cal{W}$. 
 For a $t$-tuple $x=(x_1,\ldots,x_t) \in B^t$ let $x^{-i}$ denote the $(t-1)$-tuple obtained from $x$ by deleting the $i$'th coordinate. The winning set for  $(f_1,\ldots,f_t)$ is the set of all $x=(x_1,\ldots,x_t)$ such that for all $i$ it holds that $x_i \in f_i(x^ {-i})$. Let ${\cal{S}}^{(t)}$ be the set of all strategies for the $t$ player game, and ${\cal{W}}^{(t)}$ be the set of all winning sets (note that ${\cal{W}}^{(1)}={\cal{W}}$).
 
 An alternative but equivalent way of defining strategies and winning sets for the $t$-player game is the following, which will prove useful for us. (There is a canonical isomorphism between these two definitions, stemming from the inductive definition below). The idea in the following definition is that the $t$'th player, as usual, makes her decision after viewing the points associated with the other $t-1$ players, whereas they observe her point, and as a result choose a strategy for the $(t-1)$-player game that involves only them. Formally: we proceed to define by induction. For $t=1$ we use the previous definition. For $t>1$ we view $B^t$ as  $B^{t-1} \times B$, and let $X_1=B^{t-1}$, and $X_2=B$.  We define a strategy and a winning set for $X_1 \times X_2$.
\\ A strategy for $X_1 \times X_2$ is a pair of functions $f_1,f_2$, with $f_1: X_2 \rightarrow$ $\cal{S}$$^{(t-1)}$, and
$f_2: X_1 \rightarrow$ $\cal{S}$$^{(1)}$.  The winning set for $(f_1,f_2)$ is the set of all points $x_1,x_2$ such that $x_1$ belongs to the winning set of the strategy $f_1(x_2)$ and $x_2$ belongs to the winning set of the strategy $f_2(x_1)$.

In this paper we will concentrate on $B=\{0,1\}^n$ with various choices for $\cal{W}$$^{(1)}$. We will use the uniform measure on $B^t$, which we denote by $\mu$, and define 
$$
p(t,n) := \max_{W \in {\cal{W}}^{(t)}} \mu(W),\ \ \  p(t):= \lim_{n \rightarrow \infty} p(t,n).
$$
The three choices of $\cal{W}$$^{(1)}$ that will interest us are:
\begin{itemize}
\item Let $\cal{W}$$_{dict}$ be the set of dictators, i.e. the set of all $W_i = \{x \in B : x_i=1\}$. This is the basis for the hats game. We will henceforth use $p_{dict}(t,n)$ and $p_{dict}(t)$ for $p(t,n)$ and $p(t)$, the success probabilities in this setting.
 \item Let $\cal{W}$$_{intersecting}$ be the set of all intersecting families in $\{0,1\}^n$, i.e, the set of all $W \subset \{0,1\}^n$ such that if $x,y \in W$ then there exists a coordinate $i$ such that $x_i=y_i=1$. We will use 
 $p_{intersecting}(t,n)$ and $p_{intersecting}(t)$ for the success probabilities in this case.
 \item Let $\cal{W}$ be the set of all balanced monotone families in $\{0,1\}^n$, i.e, all $W$ containing precisely half the points in $\{0,1\}^n$ such that if $x \in W$ and $y_i \geq x_i$ for all $i$ then $y \in W$.  
 We will use 
 $p_{monotone}(t,n)$ and $p_{monotone}(t)$ for the success probabilities in this case.
 \end{itemize}
 Note that every dictatorship is an intersecting family, and every maximal intersecting family is a balanced monotone family, so 
 $$
 p_{monotone}(t,n)\geq  p_{intersecting}(t,n) \geq  p_{dict}(t,n)
 $$
 and
  $$
 p_{monotone}(t)\geq  p_{intersecting}(t) \geq  p_{dict}(t).
 $$
 Thus the following two conjectures are progressively stronger than Conjecture \ref{Main}
 \begin{conjecture}\label{Kneser}
$    p_{intersecting}(t)$ tends  to 0 as $t$ grows.
\end{conjecture}
\begin{conjecture}\label{Monotone}
$    p_{monotone}(t)$ tends  to 0 as $t$ grows.
\end{conjecture}
{\bf Remark}: In all these cases, changing the order of the limits in the value of $p(t,n)$ leads to an uninteresting question, as for every fixed $n$ $p(t,n)$ tends to 0 as $t$ grows.
 \subsection{Winning sets as independent sets in Hamming products of graphs}
Having described the general setting we would like to point out that Conjecture \ref{Kneser} is actually a statement in graph theory. To that end, here are some definitions. 
\dfn
The Kneser graph $K(n)$ is a graph on vertex set $\{0,1\}^n$, with an edge between $x$ and $y$ if $x$ and $y$ have disjoint support, i.e. there is no $i$ for which $x_i=y_i=1$.
\edfn
\dfn
Let $G_1,G_2,\ldots G_t$ be graphs. The Hamming product 
$G_1 \square G_2 \square \ldots \square G_t$ has vertex set $V(G_1) \times V(G_2) \times \ldots V(G_t)$ with an edge between $x$ and $y$ if there exists an index $i \in [t]$ such that for $j \not = i$ $x_j=y_j$ and $\{x_i,y_i\} \in E(g_i)$.
So the $t$-fold product of a graph consisting of a single edge with itself is the usual graph of the Hamming cube.

To stress the relationship between the graph product and the strategies and winning sets in the $t$-player game we add also the inductive definition;
The Hamming product of graphs $G$ and $H$ has vertex set $V(G) \times V(H)$, and an edge between
$(x,v)$ and $(y,u)$ if either $x=y$ and $\{v,u\}$ is an edge in $H$, or $v=u$ and $\{x,y\}$ is an edge in $G$.
We denote it by $G \square H$. There is a canonical isomorphism between $G \square (H \square M)$ and 
$(G \square H) \square M$, so we will treat this product as an associative relation, and write $G^{\square t}$ to denote the $t$-fold Hamming product of $G$ with itself.
\edfn
\dfn
Let $\alpha(G)$ be the size of the largest independent set in $G$, and $\bar{\alpha}(G) := \frac{\alpha(G)}{|V(G)|}.$
\edfn
Note that an independent set in $G^{\square t}$ is a subset of $(V(G))^t$ such that its intersection with every 1 dimensional fiber of $(V(G))^t$ (a set of vertices resulting from fixing $t-1$ of the coordinates) is an independent set in $G$, (just as a winning set in a product space is a set whose restriction to every 1-dimensional fiber is a 1-dimensional winning set).  Also note that an independent set in $K(n)$ is an intersecting family. Thus,
\\ {\bf Observation}: 
$$
p_{intersecting}(t,n) = \bar{\alpha}(K(n)^{\square t} ).
$$
So, we may restate Conjecture \ref{Kneser} as
\begin{conjecture}\label{Kneser2}
$$
\lim_{t \rightarrow \infty} \lim_{n \rightarrow \infty} \bar{\alpha}(K(n)^{\square t} ) = 0 .
$$
\end{conjecture}
\subsection{Relating the maximal winning set in $B^{t+1}$ to the maximal winning set in a random subset of $B^t$}
 We now return to the general setting of a game on $B^t$ and $B^{t+1}$ and proceed to express $p(t+1,n)$ as the expected measure of the largest intersection of a winning set and a random set in $B^t$. 
 
 Consider the game on $B^{t+1} = B^t \times B$ and a strategy $(f_1,f_2)$, with $f_1:B \rightarrow {\cal{S}}^{(t)}$ and
  $f_2: B^t \rightarrow {\cal{S}}^{(1)}$. These two functions induce (for $i=1,2$) functions  
  $g_1:B \rightarrow {\cal{W}}^{(t)}$ and $g_2:B^t \rightarrow {\cal{W}}^{(1)}$, simply by letting $g_i(x)$ be the winning set of $f_i(x)$.

  We claim that for a given $f_2$ it is simple to describe an optimal choice of $f_1$.
  Let $\cal{W}^{(1)}$ consist of $r$ sets: 
  $ {\cal{W}}^{(1)} = \{W^{(1)}_i\}_{i=1}^r$, and
  first, note that $f_2$ (which defines $g_2$) induces a partition of $B^t$ into $V_1,\ldots V_r$, where
   $V_i:= g_2^{-1}(W^{(1)}_i)$.
  Secondly, note that a random uniform choice of $y \in B$ induces a random subset $R_y \subseteq [r]$ according to the winning sets that $y$ belongs to, i.e.
  $$
  R_y := \{ i: y \in W^{(1)}_i\}.
  $$
  Now, given $f_2$, and a fixed $x_2 \in B$, how best to define $f_1(x_2)$? Well, observe that a necessary condition for $(x_1,x_2)$ to be contained in a winning set is for $x_2 \in g_2(x_1)$ which means that $x_1 \in V_i$ for some $i \in R_{x_{2}}$. Therefore, the best choice for $f_1(x_2)$ is such that  the winning set $g_1(x_2)$ is  the winning set $W^{(t)}_j$ that maximizes the probability that a random choice of $x_1 \in B^t$   lands in  $W^{(t)}_j \cap \bigcup_{i \in R_{x_{2}}} V_i$. This implies
 \begin{lemma}\label{p(t+1)}
  $$
 p(t+1,n)= \max_{B^t=\cup_{i=1}^r V_i } E_{x_2 \in B}\left[\max_{W \in {\cal{W}}^{(t)} } 
 \mu(W \cap \bigcup_{i \in R_{x_{2}}} V_i)\right].
 $$
  \end{lemma}
 Here the first maximum is over all partitions of $B^t$, (each corresponding to a choice of $f_2$), and the second maximum represents the success probability for the optimal choice of $f_1(x_2)$, given $x_2 \in B$.

 \subsection{A special case: maximal independent sets in random subsets of hamming powers of the Kneser graph}
 We can use Lemma \ref{p(t+1)} to find an upper bound for  $p_{intersecting}(t+1,n)$ (and hence for $p_{dict}(t+1,n)$)
 in graph theoretic terms. Let ${\cal{W}}={\cal{W}}_{intersecting}=\{W_i\}_{i=1}^r$ be the family of maximal independent sets in $K(n)$, or, in other words, the family of maximal intersecting sets in $\{0,1\}^n$. A choice of a random vertex $v$ 
 in $K(n)$ induces a choice of a random set $R_v \subseteq [r]$, consisting of all indices $i$ such that $v$ belongs to $W_i$,
 $$
 R_v:=\{i : v \in W_i\}.
 $$
 Denote the distribution over $R_v$ by $\cal{D}$.
  We have, then, the following corollary of Lemma \ref{p(t+1)}.
 \beqn\label{1}
 p_{intersecting}(t+1,n)  \leq 
  \max_{ (K(n))^{\square t}=\cup_{i=1}^r V_i} 
  E_{R \sim D}\left[\max_{W \in {\cal{W}}} \mu(W \cap \bigcup_{i \in R} V_i)\right],
 \eeqn
 where we allow any partition $ (K(n))^{\square t}=\cup_{i=1}^r V_i$ into disjoint sets, some of which may be empty.

\remove{
 We note that in the proof of lemma \ref{blockers} we use the fact that the measure of every dictatorship is precisely 1/2, and that they are positively correlated. 
 Recall that every maximal intersecting family contains precisely half of the sets, (precisely one of every pair of complementary sets), so each $W \in {\cal{W}}$ has measure 1/2, so choosing $v$ uniformly at random, the marginal probability of each $i$ belonging to $R_v$ is precisely 1/2.
 Due to the positive correlation of increasing events (see, e.g., Harris' inequality  \cite{Harris}) these events are non-negatively correlated, i.e., for every $i \in [r], J \subseteq [r]$
 $$
 Pr[i \in R_v | J \subseteq R_v] \geq 1/2.
 $$
 Let ${\cal{D}}_r$ denote the set of all such distributions over $[r]$.
 We have, then, the following corollary of Lemma \ref{p(t+1)}.
 \beqn\label{1}
 p_{intersecting}(t+1,n)  \leq 
  \max_{D \in {\cal{D}}_r, (K(n))^{\square t}=\cup_{i=1}^r V_i} 
  E_{R \sim D}\left[\max_{W \in {\cal{W}}} \mu(W \cap \bigcup_{i \in R} V_i)\right],
 \eeqn
 }
 
 {\bf Remarks:}
 \begin{itemize}
 \item Equation (\ref{1}) bounds the size of the maximal independent set in the $(t+1)$'th Hamming power of the Kneser graph in terms of the maximal independent set contained in a random subset of the vertices of the $t$'th power. We will expand below on this theme, and raise some conjectures regarding this setting in general graphs.
 \item Recalling that    $p_{intersecting}(t,n) \geq p_{dictator}(t,n)$ makes this approach relevant to solving the hats problem.
 \item To solve the more general problem, using the technique of blockers, that we use for proving lemma \ref{blockers}, we would need to allow more general FKG-type distributions over $R$, where the marginal of every coordinate is Bernoulli(1/2), and there is a non-negative correlation between all sets of coordinates. 
 \end{itemize}

 \subsection{Maximal independent sets in random subgraphs}
 We would like to make a general conjecture regarding independent sets in random subgraphs, that if true, using (\ref{1}), would imply that $ \bar{\alpha}(K(n)^{\square t} )$ tends to 0 as $t$ grows, and thus also prove Levine's conjecture regarding the hats problem, Conjecture \ref{Main}.
 \\ First let us recall some definitions, and make some new ones. 
 \\ For any graph $G$ let $\mu$ denote the uniform measure on $V(G)$.
 \\ Let ${\cal{I}}(G)$ be the family of all independent sets in $G$.
 \\ Let $\bar{\alpha}(G)= \max_{I \in {\cal{I}}(G)} \mu(I)$
 \\ Let ${\cal{D}} = \cup_r {\cal{D}}_r$ denote all distributions on subsets $R$ of some finite set $[r]$, such that every $i \in [r]$ belongs to $R$ independently with probability 1/2, and all these 
 events are positively correlated. 
 \\ Let $\alpha^*(G) = \max_{r, D \in {\cal{D}}_r, V(G)=\bigcup_{i=1}^r V_i}E_{R \sim D} [ \max_{I \in {\cal{I}}(G)} \mu(I \cap (\cup_{i \in R} V_i))]$.
  \\ Let $\alpha^{**}(G) = E_W[ \max_{I \in {\cal{I}}(G)} \mu(I \cap W)]$, where $W$ is chosen uniformly over all subsets of $V(G)$.
  \\ Let $\epsilon^*(\alpha) = \inf_{G: \bar{\alpha}(G)\geq \alpha} \{ \bar{\alpha}(G) - \alpha^*(G)\}$  
  \\ Let $\epsilon^{**}(\alpha) = \inf_{G: \bar{\alpha}(G) \geq \alpha} \{ \bar{\alpha}(G) - \alpha^{**}(G)\}$
 
\begin{conjecture}\label{randomindset}
$\epsilon^*(\alpha)>0$ for all $\alpha \in (0,1/2)$.
 \end{conjecture}
This conjecture would imply Conjecture \ref{Kneser2} (or, equivalently,Conjecture \ref{Kneser})  as follows. Assume, by way of contradiction, that
 $ \bar{\alpha}(K(n)^{\square t} )$ does not tend to 0. 
 First note that the success probability is clearly non-increasing with $t$, as the first 
 $t-1$ players can try to make their guesses whilst ignoring the hats of player $t$.  since it is also bounded from below, it must tend from above to a limit, say $\alpha$. For large enough $t$ we would have 
$$,
  \alpha \leq \bar{\alpha}(K(n)^{\square t} ) < \alpha + \epsilon^*(\alpha)
$$
and hence, using (\ref{1}), and Conjecture \ref{randomindset}

$$
 \bar{\alpha}(K(n)^{\square t+1} ) \leq \alpha^*(K(n)^{\square t}) =
 $$
$$
   =\bar{\alpha}(K(n)^{\square t} ) - \epsilon^*( \bar{\alpha}(K(n)^{\square t} )) <
 \alpha + \epsilon^*(\alpha) - \epsilon^*( \bar{\alpha}(K(n)^{\square t} )) \leq \alpha,
$$
(because $\epsilon(\alpha)$ is non-decreasing. )
So $ \bar{\alpha}(K(n)^{\square t+1}) < \alpha$, contradiction.

We do not know that Conjecture \ref{randomindset} is true even in the special cases where the distribution of $I$ is simply binomial (i.e. all $i$ belong to $I$ independently with probability 1/2), or in the even more restricted case where  
each $V_i$ consists of a single vertex. Let us state this last case as a separate conjecture, as it is seems to be a fundamental problem in the study of independent sets in graphs.
\begin{conjecture}\label{alphastar}
$\epsilon^{**}(\alpha)>0$ for all $\alpha \in (0,1/2)$.
In other words:
\\ There exists a monotone non-decreasing function $\epsilon^{**}: (0,1/2) \rightarrow (0,1/2)$ such that the following holds (where the point is that $\epsilon^{**}>0$).
If $G$ is a graph on $n$ vertices with maximum independent set of size $\alpha n$, $W$ is a binomial random subset of $V(G)$, and $I_W$ is the maximal independent set contained in $W$, then 
$$
\alpha^{**}(G)= E_{W}[|I_W|/n] \leq \alpha - \epsilon^{**}(\alpha).
$$
\end{conjecture} 
So far, as stated in the following two theorems, we are able to prove this conjecture for $\alpha > 1/4$, and in the special case of regular graphs, for $\alpha >1/8$.
\begin{theorem}
\label{t16}
Let $G=(V,E)$ be a graph with $n$ vertices and independence
number $\alpha(G) = (1/4+\tau)n$, where $\tau>0$ satisfies
$\tau<1/4$. Then 
$\alpha^{**}(G) \leq 1/4+\tau-\tau^2/3$.
\end{theorem}

\begin{theorem}
\label{p17}
For any $\tau>0$  there is some $g(\tau)>0$ so that the following
holds.
Let $G=(V,E)$ be a regular graph with $n$ vertices and independence
number $\alpha(G) = (1/8+\tau)n$.
Then 
$\alpha^{**}(G) \leq 1/8+\tau-g(\tau)$.
\end{theorem}

The best upper bound that we know for $\epsilon^{**}$ comes from $G(n,p)$, with $n = \Theta(\alpha \log(1/\alpha))$, 
and a careful choice of a constant $p$ giving
$$
\epsilon^{**}(\alpha)< \alpha 2^{-\Omega(\frac{1}{\alpha})}.
$$ 
The following subsection contains the proof of theorems \ref{t16} and \ref{p17}.

\subsection{Proofs of theorems \ref{t16} and \ref{p17}}
 All logarithms in what follows are in
base $2$, unless otherwise specified. To simplify the presentation
we omit all floor and ceiling signs whenever these are not crucial.

\noindent
We will be using an old result of Hajnal \cite{Ha} (see also \cite{Ra}), relevant for us in the case where there exist independent sets containing at least half of the vertices in the graph.
\begin{proposition}\cite{Ha}
\label{Hajnal}
For every graph $G$ the cardinality of the intersection of all 
maximum independent sets plus the cardinality of the union of all
these sets is at least $2 \alpha(G)$. 
\\  Consequently, if $\alpha(G) = \alpha n$
where $\alpha>1/2$ and $n$ is the number of vertices of $G$,
this implies that there is a set of at least $(2\alpha
-1)n$  vertices contained in all maximum independent sets. 
\end{proposition}

Using Hajnal's result we next describe the proof of Theorem
\ref{t16}.
\vspace{0.2cm}

\noindent
{\bf Proof of Theorem \ref{t16}:}\,  
Without loss of generality we may assume that
$n$ is arbitrarily large, as we can replace $G$ by a union
of many vertex disjoint copies of itself and use linearity of
expectation. Assuming $n$ is large, almost every random subset of
vertices is of cardinality $(1/2+o(1))n$, hence it suffices to
show that for almost every set $W$ of $m=(1/2+o(1))n$ vertices,
the independence number of the induced subgraph of $G$ on $W$ is
smaller than $(1/4+\eps-\eps^2/2)n$. Construct the random set $W$
of size $m$ by  removing from $G$ vertices, one by one. Starting
with $V=V_0$, let $V_{i+1}$ be the set obtained from $V_i$ by 
removing a uniform random vertex of $V_i$. The set $W$ is thus
$V_{n-m}$. Let $G_i$ be the induced subgraph of $G$ on $V_i$.
Call a step $i$, $1 \leq i \leq
n-m$ of the random process above 
successful if either the independence number of 
$G_{i-1}$ is already smaller than $(1/4+\eps-\eps^2/2)n$
(note that in this case this will surely be hold in
the final graph $G_{n-m}$), or the independence number
of $G_{i}$ is strictly  smaller than that of $G_{i-1}$.
Put $i_0=(1/2-\eps)n$ and consider the graph 
$G_{i_0}$. For any
$i >i_0$, the number of vertices of $G_{i-1}$ is at most
$(1/2+\eps)n$. If its independence number is  smaller than
$(1/4+\eps-\eps^2/2)n$ then, by definition, step number $i$ is 
successful. Otherwise, by the result of Hajnal mentioned above,
the number of vertices of $G_{i-1}$ that lie in all the maximum
independent sets in it is at least $(\eps-\eps^2)n$. Since $\eps <1/4$
this is a fraction of at least  $\eps$ of the vertices of
$G_{i-1}$. Therefore, in this case,
the probability that the next chosen vertex lies in all maximum
independent sets of $G_{i-1}$ is at least $\eps$. 
We have thus shown that for every $i$ satisfying $i_0<i  \leq n-m$
the probability that
step number $i$ is successful is at least $\epsilon$. Therefore,
the probability that there are at least $\eps^2 n/2$ successful
steps during the $n-m-i_0=(\eps-o(1))n$ steps starting with
$G_{i_0}$ until we reach $G_{n-m}$  is at least the probability
that a binomial random variable with parameters
$(\eps-o(1))n$ and $\eps$ is at least $\eps^2 n/2$. This
probability is $1-o(1)$ for any fixed positive $\eps$ as 
$n$ tends to infinity. Since having that many successful steps 
ensures that the independence number of the induced subgraph
of $G$ on $W$ is at most $(1/4+\eps-\eps^2/2)n$, this completes
the proof.  \hfill $\Box$

\subsubsection{Regular graphs}
In  the proofs of Theorem \ref{p17} and later Proposition \ref{p14}  we apply
the following early version of the  container theorem of \cite{BMS} and
\cite{ST}.
\begin{theorem}[c.f. \cite{AS}, Theorem 1.6.1] 
\label{t21}
Let $G=(V,E)$ be a $d$-regular graph on $n$ vertices  and let
$\delta>0$ be a positive real. Then there is a collection
$\CC$ of subsets of $V$ of cardinality
$$
|\CC| \leq \sum_{i \leq n/\delta d} { n \choose i} 
$$
so that each $C \in \CC$ is of size at most
$\frac{n}{\delta d} +\frac{n}{2-\delta}$ and every independent set
in $G$ is fully contained in a member $C \in \CC$.
\end{theorem}

\vspace{0.2cm}

\noindent
{\bf Proof of Proposition \ref{p17}:}\, 
Let $G=(V,E)$ be as in the proposition, where $|V|=n$.
As before we may assume without loss of generality that
$n$ is sufficiently large as a function of $\eps$.
Without trying to optimize the function $g(\eps)$, let $d$ denote
the degree of regularity of $G$. 
Note that $G$ contains  a set $S$ of  at least $n/(d^2+1)$ vertices 
no two of which are adjacent or have a common neighbor.
Let $W$ be a uniform random set of vertices of $G$. If the complement
of $W$  fully contains the closed neighborhoods of $s$ vertices
of $S$, then  the independence number of the induced subgraph of
$G$ on $W$ is at most $\alpha(G)-s$. The random variable counting
the above number $s$ is a Binomial random variable with
expectation $|S|/2^{d+1} \geq \frac{n}{(d^2+1)2^{d+1}}$. 
Thus if, say, $d$ is at most $50/\eps^4$ we get that the expected
independence number of an induced subgraph of $G$ on a uniform
random set of vertices is at most 
$$
\alpha(G) -\frac{n}{(d^2+1)2^{d+1} }
$$
supplying a lower bound (of the form $2^{-\Theta(\eps^{-4})}$) 
for $g(\eps)$. We thus may and will assume that 
$d \geq \frac{50}{\eps^4}$. By Theorem \ref{t21}
with $\delta=\eps$ there is a collection $\CC$ of subsets of 
$V$, satisfying
$$
|\CC| \leq \sum_{i \leq \eps^3n/50} {n \choose i} 
\leq 2^{H(\eps^3/50) n},
$$ 
where $H(x)=-x \log x -(1-x) \log (1-x)$ is the binary entropy
function. Each member $C$ of $\CC$ is of size at most
$$
\frac{n}{\delta d} + \frac{n}{2-\delta} 
\leq \frac{\eps^3 n}{50}+ \frac{n}{2-\eps} < (\frac{1}{2}+\eps)n
$$
and every independent set of $G$ is contained in a member $C \in
\CC$.

As in the proof of Theorem \ref{t16} we can generate a random
subset $W$ of $V$ by omitting vertices one by one, starting with
$V$. Since $n$ is large almost all sets $W$ are of size $n/2+o(n)$.
Moreover, for almost all of them the size of $W \cap C$ deviates
from $|C|/2$ by at most, say, $\frac{\eps}{100} n$
for all $C \in \CC$, provided $\eps$ is sufficiently small. 
It suffices to show that with high probability the independence
number of the induced subgraph on $W$ is
at most, say, $\alpha(G)-2g(\eps) n$.
Since every
independent set is contained in at least one of the members
$C$ of $\CC$ it suffices to show that with high probability 
the independence number of the induced subgraph of $G$ on $W \cap
C$ is at most $\alpha(G)-2g(\eps)n$ for every $C \in \CC$. 
Fix $C \in \CC$. Without loss
of generality its size is at least $n/8$ (since otherwise it
cannot contain a large independent set at all). Recall that
$|C| \leq (1/2+\eps)n$. If the independence number of the induced
subgraph of $G$ on $C$ is smaller than $(1/4+\eps)|C|$ then so is
the independence number of the induced subgraph on $W \cap C$,
and this is smaller than $\alpha(G) -0.1 \eps n$, implying the
desired result.
Otherwise, as in the proof of Theorem \ref{t16}, in the random 
process that omits vertices of $C$ one by one to get
$W \cap C$, the number of times the independence number
drops dominates stochastically a binomial random variable 
with parameters $\frac{\eps}{2} |C|$ and $\eps$. By the standard
estimates for Binomial distributions (c.f., e.g., \cite{AS},
Theorem A.1.13), the probability this  variable is less
than half its expectation is at most
$$
e^{- \eps^2 |C|/16} \leq e^{-\eps^2 n/128}.
$$
By the union bound over all $C \in \CC$ the probability this 
happens even for a single $C \in \CC$ is at most
$$
2^{H(\eps^3/50)n} \cdot e^{-\eps^2 n /128}
$$
which, for small $\eps$, tends to $0$ as $n$ tends to infinity.
This shows that in this case  ($d \geq \frac{50}{\eps^4}$ ),
with high probability the independence number
of the induced subgraph of $G$ on $W \cap C$ is smaller than
$\alpha(G)$ by at least, say, $\eps^2  n/40$, completing the 
proof of the proposition.  \hfill $\Box$

\section{Blockers and proof of the main theorem}
\subsection{Bounding $p(t+1)$ using blockers}

We now focus on the hats game, i.e. consider $B =\{0,1\}^n$, with winning sets $W_i = \{x \in B : x_i=1\}$ for $i \in \{1,\ldots,n\}$.
Let $\mu$ denote the uniform measure on $B$, and by abuse of notation, also on $B^t$.
 Call a subset of $B^t$ a winning set, if it is the winning set of any strategy for the corresponding game. Write $p(t,n)$ and $p(t)$ for short for $p_{dict}(t,n)$ and $p_{dict}(t).$
\dfn
A blocker $A \subset B^t$ is a set of points that intersects every winning set. 
\edfn

\begin{lemma}\label{blockers}
If there exist disjoint blockers $A_1,\ldots, A_r \subset B^t$, such that
\begin{enumerate}
\item $|A_i| =k$ for all $i$.
\item $\mu(\bigcup A_i) = \beta.$
\end{enumerate}
Then $p(t+1) \leq p(t)-2^{-k}\beta/k$.
\end{lemma}
{\bf Proof} : By Lemma \ref{p(t+1)} we know that $p(t+1,n)$ is bounded by the expectation of the maximal intersection of any winning set in $B^t$ with $V$, a random binomial union of subsets of $B^t$. 
 Now, if a blocker $A$ is disjoint from $V$ this means that in {\em every} winning set $W$  there is at least one point from $A$ (hence not in $V$), 1 point being $1/k$ of the measure of $A$. If the union of all the blockers missed by $V$ has measure $\tau$ this means every winning set contains a set of measure at least $\tau/k$ that's disjoint from $V$. Now, since every blocker is missed with probability at least $2^{-k}$, the expected proportion of missed blockers is $2^{-k}$ , contributing $2^{-k}\beta$ to the expected measure of the union of missed blockers, which contributes 
 at least $2^{-k}\beta/k$ to the expected measure-loss of the maximum that defines $p(t+1,n)$. Taking the limit with respect to $n$ gives the required result.

\subsection{Constructing blockers for the hats game}
In this subsection we consider the hats game, and construct, for every $t$, by induction on $t$, a set of blockers for $B^t$.
This, together with Lemma \ref{blockers}, will prove the main claim of this paper, Theorem \ref{thm:main}
\begin{lemma}
 Let $k(1)=2$, $\beta(1)=1$, and for $d\ge1$ define $k(d+1):= k(d) {{2k(d)}\choose{k(d)}}$. (So $k$ grows as a tower function of $d$), and $\beta(d+1) = \beta(d) (2 {{2k(d)}\choose{k(d)}} )^{-1}.$
 Then, for every $d$ there exist a family of blockers $A_1,\ldots,A_r \subset B^d$ with
 \begin{enumerate}
\item $|A_i| =k(d)$ for all $i$.
\item $\mu(\bigcup A_i) \ge \beta(d).$
\end{enumerate}
 \end{lemma}

\begin{corollary}
For all $d > 1$ 
$$
p(d+1) \le p(d)-2^{-k(d)}\beta(d)/k(d) < p(d).
$$
\end{corollary}

{\bf Proof of Lemma}
We will build the family of blockers for $B^d$ inductively.
For $d=1$ the set of blockers for $B$ is the set of all pairs $\{x,\bar{x}\}$. Every dictator must contain precisely one element from each pair. Now, assume we have a family of blockers of size $k(d)$ for $B^d$ as desired.
Let $\ell : = {{2k(d)}\choose{k(d)}}$. We will choose randomly (in a manner to be described below) a series of disjoint unordered $\ell$-tuples
$Y^{(j)}=\{y^{(j)}_1,y^{(j)}_2,\ldots,y^{(j)}_\ell\}$, with $y^{(j)}_i \in B$, for $j = 1,2,\ldots$ until the measure of the union of these $\ell$-tuples in $B$ exceeds
$\frac{1}{2\ell}$. The new blockers for $B^{d+1}$ will be all cartesian products of the form $b \times Y^{(j)}$ where $b$ is one of the blockers we designed for $B^d$. The claims regarding the size of the blockers and the measure of their union are immediate. We must check two things. First, that the product  $b \times Y^{(j)} $ is indeed a blocker for $B^{d+1}$, secondly, that one can choose the desired number of disjoint $\ell$-tuples.
\\ To this end, let us describe how the $\ell$-tuples are formed. We take a random partition of the $n$ coordinates of $B$ into $2k(d)$ sets $S_1,\ldots,S_{2k(d)}$, uniformly over all such (ordered) partitions. For every $I:= \{ i_1,\ldots,i_{k(d)} \}\subset [2k(d)]$, let $y_I \in B$ be the vector whose 1-support is precisely $S_{i_1}\cup S_{i_2}\cup\ldots\cup S_{i_{k(d)}}$. This defines $\ell$ different vectors corresponding to the specific partition. We proceed to choose such  $\ell$-tuples sequentially at random, and discard any $\ell$-tuple that is not disjoint from all its predecessors. Note that the marginal distribution of every $Y_I$ is uniform, hence if the union of all predecessors of $Y^{(j)}$ has measure $\epsilon$ then with probability at least $1-\ell \epsilon$ it will be disjoint from its predecessors, so we may, as claimed, continue until the measure of the union of all $\ell$-tuples is $\frac{1}{\ell}(1-o(1))$.
Finally, we prove that $b \times Y^{(j)} $ is a blocker for $B^d \times B$, where $Y^{(j)}$ is an $\ell$-tuple corresponding to some partition, and $b = \{ x_1,\ldots, x_{k(d)} \}$ is a blocker for $B^d$. 
Let $(f_1,f_2)$ be a strategy for $B^d \times B$, and $(g_1, g_2)$ the corresponding functions such that $g_i(z)$ is the winning set for the strategy $f_i(z)$. For every $x \in B^d$, $g_2$ picks a dictatorship $g_2(x)=W_i$,
so define $t_1,\ldots,t_{k(d)}$ by $g(x_i)=W_{t_{i}}$.
These are, respectively, the dictators that the second player guesses when she sees one of the $x$'s from $b$ on the first player's forehead. For $i \in \{1,\ldots, k(d)\}$ let $S_{r_i}$ be the part of the partition of $[n]$ (used to define $Y^{(j)}$) that contains $t_i$, and let $I = \{i_{r_1},\ldots,i_{r_{k(d)}}\}$ (or an arbitrary  set of size $k(d)$ containing it if the elements in it are not distinct).  So the $\ell$-tuple $Y^{(j)}$ contains a vector $y_I$ for which all the coordinates $t_1,\ldots,t_{k(d)}$ are equal to 1, meaning that if the second player has $y_I$ on her forehead she will make a correct guess if the first player has any of the $x$'s in $b$ on her forehead, i.e. 
$$
y_I \in \cap_{i=1}^{k(d)} W_{t_i} =\cap_{x_i \in b} g_2(x_i)
$$
Now, let $g_1(y_I)$ be the corresponding winning set in $B^d$ that the first player guesses when seeing $y_I$. By the fact that $b$ is a blocker there exists $x_i \in b$ that belongs to $g_1(y_I)$, (and, as mentioned, $y_I$ belongs to the winning set $g_2(x_i)$), so the pair $(x_i,y_I)$ belongs to the winning set of $(f_1,f_2)$ - i.e. we have proven that $b\times Y^{(j)}$ intersects every winning set for $B^d \times B.$ \qed  

\section{Blockers (hitting sets) in graphs}\label{graphs}
In light of the partial success in the previous section, a tempting approach to conjecture \ref{randomindset} is to try and prove the existence of a family of disjoint blockers (to the family of maximal independent sets of almost maximum size) in any graph $G$, where their size and the measure of their union is a function of $\bar{\alpha}(G)$ . However, this conjecture, which constitutes a radical strengthening of a conjecture of Bollob\'as, Erd\H{o}s and Tuza, (see below), is too good to be true, as shown by the two (families of) examples in this section, where every blocker has size which tends to infinity as the size of the graph in the example grows. In the example of Theorem \ref{t12} the maximum-size independent sets are of size $n/4$, and the smallest set intersecting all maximal independent sets is of size $\Theta(\sqrt{n})$. In the example of theorem \ref{t13} the maximum-size independent sets are of size 
$n (1/2 -o(1))$ and the smallest set intersecting all maximal independent sets has size $\Theta (\log(n))$.

\vspace{0.2cm}

For a graph $G=(V,E)$ let $h(G)$ denote the minimum cardinality of
a set of vertices that intersects every maximum independent set of
$G$. Bollob\'as, Erd\H{o}s and Tuza (see \cite{Er}, page 224, or
\cite{CG}, page 52) raised the following conjecture.
\begin{conjecture}[\cite{Er}, \cite{CG}]
\label{cbet}
For any positive $\alpha$,
if the size $\alpha(G)$ of a
maximum independent set in an $n$-vertex graph $G$ 
is at least $\alpha n$, then $h(G)=o(n)$.
\end{conjecture}
Our examples, related to this conjecture, are as follows. As we will see, in the case of regular graphs, the first example gives an almost tight bound for its parameter range.

\begin{theorem}
\label{t12}
For every positive integer $k$ there is a graph $G=G_k$ with
$n=2k(2k-1)$ vertices, independence number $\alpha(G)=k^2 (> n/4)$,
and $h(G)=k+1 (> \sqrt {n}/2).$
\end{theorem}
\begin{theorem}
\label{t13}
For any positive integers $m$ and $t$, where $m$ is even and $4t^2
\leq m$,
there is a graph $G=G_{m,t}$ on $n=2^m$ vertices with independence number 
$\alpha(G)=\sum_{i=0}^{m/2-t} { m \choose i}$ and
$h(G)=\Theta(t^2)$.
\end{theorem}
Theorems \ref{t12} and \ref{t13} are proven in subsection{constructions} below.

We remark that Theorem \ref{t13} settles 
the final open problem raised by Dong and Wu in \cite{DW}.

The graphs establishing the assertion of Theorem \ref{t12} are
regular. It turns out that for regular graphs (of any degree)
the estimates in this theorem are nearly tight, as stated 
in the next proposition.
\begin{proposition}
\label{p14}
For any fixed $\eps>0$ and any regular graph $G$ with $n>n_0(\eps)$ 
vertices
satisfying $\alpha(G)  \geq (1/4+\eps)n$, the parameter $h(G)$
satisfies $h(G) < (1/\eps) \sqrt {n \log n}+1$.
\end{proposition}

\noindent
{\bf Proof of Proposition \ref{p14}:}\, 
Let $G=(V,E)$ be a $d$-regular graph on $n$ 
vertices with independence number at least $(\frac{1}{4}+\eps)n$,
and assume that
$n$ is sufficiently large as a function of $\eps$.
The closed neighborhood
of any vertex of $G$ intersects  every maximum independent set
of $G$, implying that $h(G) \leq d+1$.  If 
$d \leq (1/\eps) \sqrt {n \log n}$ this
implies the desired result, hence we may and will assume that
$d$ is larger. By Theorem \ref{t21} with $\delta=\epsilon$ 
there is a collection $\CC$ of at most
$$
\sum_{i \leq \sqrt n/ \sqrt {\log n} }
{ n \choose i} \leq 2^{\sqrt {n \log n} }
$$
subsets of $V$, each of size at most
$$
\frac{n}{\sqrt {n  \log n}} +\frac{n}{2-\eps}<
(\frac{1}{2}+\eps)n
$$ 
so that every independent set of $G$ is fully
contained in one of them. 

Let $X$ be a random set of $\frac{1}{\eps} \sqrt {n \log n}$
vertices chosen uniformly (with repetitions) among all
vertices of $G$. Fix a container $C \in \CC$. By Hajnal's result, Proposition \ref{Hajnal},
there are at least $\eps n$ vertices contained in all 
maximum independent sets of $G$ that are contained in $C$. The
probability that $X$ does not contain any of these vertices
is at most 
$$(1-\eps)^{\frac{1}{\eps} \sqrt {n \log n}}
\leq e^{-\sqrt {n \log n}}.
$$
The desired result follows by applying the 
union bound over all $C \in \CC$.
\hfil $\Box$

\subsection{Constructions}

\noindent
{\bf Proof of Theorem \ref{t12}:}\, 
The graph $G=G_k$ is the shift graph described as follows.
Put $K=\{1,2,.., 2k\}$. The  set of 
vertices of $G_k$ is the set of all
ordered pairs $(i,j)$ with $i \neq j$ and $i,j \in K$. Thus the 
number of vertices is 
$n=2k(2k-1)$. Two vertices $(a,b)$ and $(c,d)$ are adjacent
if $b=c$ or $d=a$.
Note that the vertices can be viewed as all
directed
edges of the complete directed graph on $K$, where 
two are adjacent iff they form
a directed path of length $2$. It is easy to check 
that the maximum independent sets of this graph are of size
$k^2$. Indeed, for every partition of $K$ into two disjoint 
parts $S$ and $T$
of equal cardinality, the set of all
pairs $(s,t)$ with $s \in S, t \in T$ is a maximum independent
set, and these are all the maximum independent sets.
Any set $H$ of at most $k$ vertices of $G$ can be viewed as
$k$ directed edges of the complete graph on $K$. Let $S$
be a set of $k$ points in $K$ that does not contain the 
head of any of these $k$ directed edges, and put $T=K-S$.
Then the maximum independent set consisting of all
pairs $(s,t)$ with $s \in S, t \in T$ does not intersect $H$.
Therefore $h(G) \geq k+1$. This is tight as shown by 
a set of pairs forming a directed cycle of length $k+1$ in the
complete directed graph on $K$. \hfill $\Box$
\vspace{0.2cm}

\noindent
{\bf Proof of Theorem \ref{t13}:}\, 
Let $G=G_{m,t}$ be the graph whose vertices are all binary vectors
of length $m$, where two are adjacent iff the Hamming distance
between them exceeds $m-2t$. Note that this is the Cayley graph
of $Z_2^m$ with respect to the set of all vectors of Hamming weight
at least $m-2t+1$. This graph contains as an induced subgraph
the Kneser graph $K(m,m/2 +1 -t)$.
By an old result of Kleitman \cite{Kl},
the independence number of this graph is
exactly $\sum_{i=0}^{m/2-t} {m \choose i}$.
The maximum independent sets are the $2^m$ Hamming
balls of radius $m/2-t$ centered at the vertices of $G$. Any set of
vertices that hits all these independent sets forms a covering code
of radius $m/2-t$ in $Z_2^m$. By using known results about covering
codes  in this range of the parameters 
it is not
difficult to prove that the minimum possible size of such a set is
$\Omega(t^2)$. Indeed, viewing the vectors of the covering code
as vectors with $\{-1,1\}$ coordinates, if their number is $T$ then 
by a known result in Discrepancy Theory (see, e.g., \cite{AS}, 
Corollary 13.3.4), there is a $\{-1,1\}$ vector whose inner product 
with all members of the code is in absolute value at most
$12 \sqrt T$. If $12 \sqrt T < 2t$ this gives a vector 
whose Hamming distance from any codeword is larger than $m/2-t$,
contradicting the assumption. This shows that the size of the code
is at least $\Omega(t^2)$. 
This is tight up to the hidden constant in the
$\Omega$-notation as can be shown by a random construction of
vectors of length $\Theta(t^2)$, extending each such vector in two
complementary ways on the remaining coordinates, or by taking the
rows of a Hadamard matrix of order $\Theta(t^2)$ and their
inverses, extending them in the same way.
Note that the fact that the Kneser graph
$K(m,m/2+1 -t)$ is a subgraph of $G$ also
implies a lower bound of $2t$ for the size of the hitting set
(as the Hamming balls  of radius $m/2-t$ centered in the points of
the hitting set cover all points, providing a proper coloring of the
Kneser graph), but the bound obtained this way is weaker than the
tight $\Theta(t^2)$ bound. \hfill  $\Box$

\section{Remarks and a conjecture}
\begin{itemize}
\item
Conjecture \ref{cbet} remains open for $n$-vertex graphs
with independence number at most $n/2$ and for such regular graphs
of independence number at most $n/4$. Similarly,
Conjecture \ref{alphastar} remains open for 
$n$-vertex graphs
with independence number at most $n/4$ and for such regular graphs
with independence number at most $n/8$. Both conjectures
appear to be significantly more difficult for graphs with
independence number $\beta n$ when $\beta>0$ is a fixed  small
positive real.  
\item As shown by the two constructions in the previous section, one cannot hope to prove Conjecture \ref{alphastar} by constructing a large family of bounded-size blockers, as in the proof of theorem \ref{thm:main}. However, in these two examples, where the minimal size of a blocker (hitting set) is large, the number of maximum independent sets is very small. In the shift graph, where the blockers are of size $\Theta(\sqrt{n})$ the number of maximum independent sets is of order $2^{\Theta(\sqrt{n})}$. In the second example the number of maximum independent sets is only $n$.
This leads to the following conjecture, asserting that the family of independent sets may be partitioned into a small number of parts, where for each part the maximal independent sets of almost maximum-size may be blocked by a large family of disjoint blockers of bounded size. This would suffice to imply Conjecture \ref{alphastar}.
\begin{conjecture}
\label{bold}
For every $\alpha >0$ there exists a positive integer $B$ and $\tau>0$ and $\tau4^{-B} > \delta >0$, and $\epsilon>0$ such that the following holds for sufficiently large $n$.
Let $G$  be a graph on $n$ vertices, where the maximum independent sets are of size $\alpha n$, and let  ${\cal{I}}$ be the family of all maximal independent sets in $G$ of size at least $(\alpha-\epsilon)n$.
Then $\cal{I}$ can be partitioned into at most $2^{\delta n}$ parts, where for each part $I_j$ there exists a family of  $\tau n$ pairwise disjoint sets of size $B$ that each intersect every element in $I_j$ .
\end{conjecture}
\item
A conjecture raised by the first author more than ten years ago motivated 
by some of the results in \cite{AHLSW}
is that the chromatic number of the graph $G_{m,t}$ described
in the proof of Theorem \ref{t13}, where $4t^2 \leq m$,
is $\Theta(t^2)$. This has been mentioned in several lectures,
see, for example, \cite{Al}. By the arguments described in the
proof of Theorem \ref{t13} this chromatic number is at least
$2t$ and at most $O(t^2)$.
\end{itemize}

\section{Acknowledgments} We thank Wojtech Samotij for useful discussions, and Zichao Dong and
Zhuo Wu for telling us about \cite{DW}.


\begin{thebibliography}{99}
\bibitem{Al} N. Alon, Graph Coloring: Local and Global, Public
Lecture, Harvard, 2017, 
https://www.youtube.com/watch?v=lFD\_DeWodn8
\bibitem{AHLSW}
N. Alon, A. Hassidim, E. Lubetzky, U. Stav and A. Weinstein,
Broadcasting with side information,
Proc. of the $49^{th}$ IEEE FOCS (2008), 823-832.
\bibitem{AS}
N. Alon and J. H. Spencer, The Probabilistic Method,
Fourth Edition, Wiley, 2016, xiv+375 pp.
\bibitem{AT} N. Alon, G. Tardos, Private communication. 
\bibitem{BMS}
J. Balogh, R. Morris, and W. Samotij,
Independent sets in hypergraphs,
J. Amer. Math. Soc. 28(2015), 669--709.

\bibitem{BFGKRTVY}   J. Buhler, C. Freiling, R. Graham, J. Kariv, J. R. Roche, M. Tiefenbruck, C. Van Alten, D. Yeroshkin,  {\em On Levine's notorious hat puzzle}, arXiv:1407.4711, 2021.

\bibitem{CG}
F. Chung and R. L. Graham, Erd\H{o}s on Graphs, His Legacy of Unsolved
Problems, A K Peters, Ltd., Wellesley, MA, 1998. xiv+142 pp.
\bibitem{DW}
Z. Dong and Z. Wu,
On the stability of graph independence number, 
arXiv:2102.13306v2, 2021.
\bibitem{Er}
P. Erd\H{o}s,
Problems and results on set systems and hypergraphs, Extremal
problems for finite sets (Visegr\'ad, 1991), 217--227, Bolyai
Soc. Math.  Stud., 3, J\'anos Bolyai Math. Soc., Budapest, 1994
\bibitem{FL13} T. Friedrich and L. Levine. Fast simulation of large-scale growth models. Random Structures and  Algorithms, 42:185-213, 2013.
\bibitem{Ha}
A Hajnal, 
A theorem on k-saturated graphs, Canadian J. Math.,
17 (1965), 720--724.

\bibitem{Harris}
T. E. Harris, Lower bound for the critical probability in a certain percolation process, Math. Proc. Cambridge Phil. Soc. 56 (1960) 13--20


\bibitem{Khova} T.Khovanova. How many hats can fit on your head?, 2011. blog.tanyakhovanova.com $/2011/04$

\bibitem{Kl}
D. J. Kleitman, 
On a combinatorial conjecture of Erd\H{o}s, J. Combinatorial Theory 1
(1966), 209--214.
\bibitem{Ra}
L. Rabern, On hitting all maximum cliques with an independent set,
arXiv:0907.3705, 2009.
\bibitem{ST}
D. Saxton and A. Thomason, Hypergraph containers,
Invent. Math. 201 (2015), 925--992.

\end{thebibliography}
 \end{document}